\documentclass[12pt]{amsart}
\usepackage{amssymb,amsmath,amscd,oldgerm}
\newtheorem{Th}{\sc Theorem}[section]

\newtheorem{Prop}[Th]{\sc Proposition}
\newtheorem{Lem}[Th]{\sc Lemma}

\newtheorem{Cor}[Th]{\sc Corollary}
\newtheorem{Rem}[Th]{\sc Remark}

\DeclareMathOperator{\des}{Des}

\begin{document}

\title{\bf On the Eulerian Polynomials of Type $D$}

\author{Chak-On~Chow}
\address{Department of Mathematics\\
The Hong Kong University of Science \& Technology\\
Clear Water Bay, Kowloon, Hong Kong}
\email{\tt cchow@ust.hk}

\maketitle

\begin{abstract}
We define and study sub-Eulerian polynomials of type $D$,
which count the elements of $D_n$ with respect to the number of
descents in a refined sense.
The recurrence relations and exponential generating functions of the
sub-Eulerian polynomials are determined,
by which the solution to a problem of Brenti,
concerning the recurrence relation for the Eulerian polynomials of type $D$,
is also obtained.
\end{abstract}

\section{Introduction}

One of the classical polynomials of combinatorial significance is the
Eulerian polynomial, which enumerates elements of the symmetric group
with respect to the number of descents, and whose properties are well
studied (see, e.g., \cite{s96}).

\bigskip
The Eulerian polynomials and their $q$-generalizations
have also been defined for other Coxeter families (see, e.g., \cite{b94}
for $q$-Eulerian polynomials which interpolate between the type
$A$ and $B$, and between type $A$ and $D$, cases).
For instance, the basic properties of the Eulerian polynomials of type $B$,
analogous to their type $A$ counterparts, are known.

\bigskip
On the other hand, the type $D$ theory is not as well developed as the
type $B$ one does.
Only the type $D$ generating functions and Worpitzky identity
are known; the type $D$ recurrence relation is still missing.
Finding such a type $D$ recurrence happens to be a problem not as simple as in
the type $A$ and $B$ cases.
The difficulty lies in that, by imitating the derivations of the
corresponding type $A$ and $B$ recurrences, the number of descents changes
in a peculiar way so that some refined Eulerian polynomials are needed in order
to capture these changes. This leads to the introduction of the sub-Eulerian
polynomials, which is the idea central to the present work.

\bigskip
The organization of this paper is as follows.
In the next section, we collect the notations that will be used in the rest
of this work. In section 3, we introduce the sub-Eulerian polynomials,
which enumerate the elements of the Coxeter group of type $D$ in a refined
sense. In particular, we obtain a recurrence relation involving the type $D$
Eulerian and sub-Eulerian polynomials.
In section 4, we compute the exponential generating functions of
the sub-Eulerian polynomials. The generating functions computed enable
us to refine the recurrence relation obtained in section 3.
In the final section, we determine the partial differential equation
of minimal order which the exponential generating function of type $D$
Eulerian polynomials satisfy,
and by which a recurrence involving the type $D$ Eulerian
polynomials and its derivatives only is derived.

\section{Notations}

We collect some definitions, and notations that will be used in the rest of
this paper.
Let $S$ be a finite set. Denote by $\#S$ the cardinality of $S$.
Denote by $\mathfrak{S}_n$ the symmetric group of degree $n$,
$B_n$ the hyperoctahedral group of rank $n$,
and $D_n$ the Coxeter group of type $D$ of rank $n$.
Let $\pi=\pi_1\pi_2\cdots\pi_n\in D_n$, where $\pi(i)=\pi_i$,
for $i=1,2,\ldots,n$.
We say that $i\in [0,n-1]$ is a $D$-descent of $\pi$ if $\pi_i>\pi_{i+1}$,
where $\pi_0=-\pi_2$.
We shall drop the type designation of descents in the sequel unless
circumstances demand the contrary.
Denote by $\des(\pi)$ the descent set of $\pi$,
and $d(\pi)=\#\des(\pi)$ the number of descents of $\pi$.
Denote by $D_{n,k}$ the Eulerian number of type $D$,
which is defined as the number of elements of $D_n$ with $k$ descents.
Define the Eulerian polynomial $D_n(t)$ of type $D$ by
$$
D_n(t)=\sum_{\pi\in D_n}t^{d(\pi)}=\sum_{k=0}^n D_{n,k}t^k.
$$

\section{Sub-Eulerian Polynomials}

The goal of the present work is to study the recurrence relation satisfied
by $D_{n,k}$.
Towards this end, we have the symmetric group $\mathfrak{S}_n$,
and the hyperoctahedral group $B_n$ as our guiding examples.
Recall that $\mathfrak{S}_n$ (resp., $B_n$) can be constructed by inserting
$n$ (resp., $\pm n$) to the elements of $\mathfrak{S}_{n-1}$
(resp., $B_{n-1}$).
By studying how the descent number changes,
recurrence relations for the Eulerian numbers of type $A$ and $B$ are obtained.

\bigskip
In the case of $D_n$, a similar construction is also possible.
Let $\sigma=\sigma_1\sigma_2\cdots\sigma_{n-1}\in B_{n-1}$.
Denote by $\hat\sigma$
the element $\bar\sigma_1\sigma_2\cdots\sigma_{n-1}$ of $B_{n-1}$.
It is clear that the map $B_{n-1}\to B_{n-1}$, $\sigma\to\hat\sigma$,
is a bijection sending $D_{n-1}$ onto $B_{n-1}\setminus D_{n-1}$.
To each element $\sigma$ of $D_{n-1}$, we can associate to $\sigma$
the pair of elements $\{\sigma,\hat\sigma\}$ of $B_{n-1}$.
This can also be realized as the orbit of $\sigma$ under the right action
of the subgroup $\langle s_0\rangle$ of $B_{n-1}$ generated by
$s_0=(1\ \bar 1)$ (in cycle notation).
(It is unfortunate that this right action of $\langle s_0\rangle$ does not
endow the collection of orbits with a group structure because
$\langle s_0\rangle$ is not normal.)
With this association in place, $D_n$ is constructible from $D_{n-1}$
by inserting $n$ to $\sigma$ and $-n$ to $\hat\sigma$.

\bigskip
For $n\geqslant 2$,
denote by ${}^1D_{n,k}$ (resp., ${}^{0,1}D_{n,k}$, ${}^{\geqslant 2}D_{n,k}$,
and ${}^{0,\geqslant 2}D_{n,k}$) the collection of elements of $D_n$ of $k$
$D$-descents, with descent at position $1$ but not at $0$
(resp., with descents at both positions $0$ and $1$,
with no descent at positions $0$ and $1$,
and with descent at position $0$ but not at $1$).
Denote by ${}^1\bar D_{n,k}$ (resp., ${}^{0,1}\bar D_{n,k}$,
${}^{\geqslant 2}\bar D_{n,k}$, and ${}^{0,\geqslant 2}\bar D_{n,k}$)
the collection of elements of $B_n\setminus D_n$ of $k$ $D$-descents,
with descent at position $1$ but not at $0$
(resp., with descents at both positions $0$ and $1$,
with no descent at positions $0$ and $1$,
and with descent at position $0$ but not at $1$).

\bigskip
Let $\sigma=\sigma_1\sigma_2\cdots\sigma_n\in D_n$ be of $k$ descents.
It is of interest to see how the pair $\{\sigma,\hat\sigma\}$ distributes
amongst ${}^1D_{n,k}$, ${}^1\bar D_{n,k}$, etc.

\bigskip
\begin{Lem}\label{lem31}
We have
$$
\begin{array}{rl}
\sigma\in{}^1D_{n,k}\iff&\hat\sigma\in{}^{0,\geqslant 2}\bar D_{n,k},\\
\sigma\in{}^{0,1}D_{n,k}\iff&\hat\sigma\in{}^{0,1}\bar D_{n,k},\\
\sigma\in{}^{\geqslant 2}D_{n,k}\iff&
\hat\sigma\in{}^{\geqslant 2}\bar D_{n,k},\\
\sigma\in{}^{0,\geqslant 2}D_{n,k}\iff&\hat\sigma\in{}^1\bar D_{n,k}.
\end{array}
$$
\end{Lem}
\begin{proof}
The map $\sigma\to\hat\sigma$ does not alter $\sigma_2\cdots\sigma_n$
so that descents beyond the position $1$ are unaffected.
It remains to figure out how descents at positions $0$ and $1$ change.
But the following calculations
$$
\begin{array}{rrrl}
\sigma\in {}^1D_{n,k}\iff&
\left.\begin{matrix}
\sigma_1+\sigma_2>0\\
\sigma_1>\sigma_2
\end{matrix}\right\}\iff&
\left.\begin{matrix}
\bar\sigma_1<\sigma_2\\
0>\bar\sigma_1+\sigma_2
\end{matrix}\right\}\iff&
\hat\sigma\in {}^{0,\geqslant 2}\bar D_{n,k},\\
\sigma\in {}^{0,1}D_{n,k}\iff&
\left.\begin{matrix}
\sigma_1+\sigma_2<0\\
\sigma_1>\sigma_2
\end{matrix}\right\}\iff&
\left.\begin{matrix}
\bar\sigma_1>\sigma_2\\
0>\bar\sigma_1+\sigma_2
\end{matrix}\right\}\iff&
\hat\sigma\in {}^{0,1}\bar D_{n,k},\\
\sigma\in {}^{\geqslant 2}D_{n,k}\iff&
\left.\begin{matrix}
\sigma_1+\sigma_2>0\\
\sigma_1<\sigma_2
\end{matrix}\right\}\iff&
\left.\begin{matrix}
\bar\sigma_1<\sigma_2\\
0<\bar\sigma_1+\sigma_2
\end{matrix}\right\}\iff&
\hat\sigma\in {}^{\geqslant 2}\bar D_{n,k},\\
\sigma\in {}^{0,\geqslant 2}D_{n,k}\iff&
\left.\begin{matrix}
\sigma_1+\sigma_2<0\\
\sigma_1<\sigma_2
\end{matrix}\right\}\iff&
\left.\begin{matrix}
\bar\sigma_1>\sigma_2\\
0<\bar\sigma_1+\sigma_2
\end{matrix}\right\}\iff&
\hat\sigma\in {}^1\bar D_{n,k}
\end{array}
$$
then show that either the descents at positions $0$ and $1$ are preserved,
or the descent at position $0$ of $\sigma$ is mapped to the descent at
position $1$ of $\hat\sigma$, and conversely.
In any case, the number of descents are preserved,
and the leading descent structure are given as in the lemma,
concluding the proof.
\end{proof}

\bigskip
We are now ready to look at the insertion process.
Let $\sigma=\sigma_1\sigma_2\cdots\sigma_{n-1}\in {}^*D_{n-1,k}$,
where ${}^*D_{n-1,k}$ is any of the subcollections of elements of $D_{n-1}$
of $k$ descents, defined above.
The letters $\pm n$ can be inserted at any of the $n$ positions,
numbered consecutively from $0$ to $n-1$, and the insertion proceeds
as follows.

\bigskip\noindent
At position $0$:

if $\sigma\in {}^1D_{n-1,k}$
($\hat\sigma\in {}^{0,\geqslant 2}D_{n-1,k}$), then
\begin{equation}\label{i1}
n\sigma_1\sigma_2\cdots\sigma_{n-1}\in {}^1D_{n,k+1},\qquad
\bar n\bar\sigma_1\sigma_2\cdots\sigma_{n-1}\in {}^{0,\geqslant 2}D_{n,k};
\end{equation}
\qquad if $\sigma\in {}^{0,1}D_{n-1,k}$
($\hat\sigma\in {}^{0,1}D_{n-1,k}$), then
\begin{equation}\label{i2}
n\sigma_1\sigma_2\cdots\sigma_{n-1}\in {}^1D_{n,k},\qquad
\bar n\bar\sigma_1\sigma_2\cdots\sigma_{n-1}\in {}^{0,\geqslant 2}D_{n,k};
\end{equation}
\qquad if $\sigma\in {}^{\geqslant 2}D_{n-1,k}$
($\hat\sigma\in {}^{\geqslant 2}D_{n-1,k}$), then
\begin{equation}\label{i3}
n\sigma_1\sigma_2\cdots\sigma_{n-1}\in {}^1D_{n,k+1},\qquad
\bar n\bar\sigma_1\sigma_2\cdots\sigma_{n-1}\in {}^{0,\geqslant 2}D_{n,k+1};
\end{equation}
\qquad if $\sigma\in {}^{0,\geqslant 2}D_{n-1,k}$
($\hat\sigma\in {}^1D_{n-1,k}$), then
\begin{equation}\label{i4}
n\sigma_1\sigma_2\cdots\sigma_{n-1}\in {}^1D_{n,k},\qquad
\bar n\bar\sigma_1\sigma_2\cdots\sigma_{n-1}\in {}^{0,\geqslant 2}D_{n,k+1}.
\end{equation}
At position $1$:

if $\sigma\in {}^1D_{n-1,k}$
($\hat\sigma\in {}^{0,\geqslant 2}D_{n-1,k}$), then
\begin{equation}\label{i5}
\sigma_1 n\sigma_2\cdots\sigma_{n-1}\in {}^{\geqslant 2}D_{n,k},\qquad
\bar\sigma_1\bar n\sigma_2\cdots\sigma_{n-1}\in {}^{0,1}D_{n,k+1};
\end{equation}
\qquad if $\sigma\in {}^{0,1}D_{n-1,k}$
($\hat\sigma\in {}^{0,1}D_{n-1,k}$), then
\begin{equation}\label{i6}
\sigma_1 n\sigma_2\cdots\sigma_{n-1}\in {}^{\geqslant 2}D_{n,k-1},\qquad
\bar\sigma_1\bar n\sigma_2\cdots\sigma_{n-1}\in {}^{0,1}D_{n,k};
\end{equation}
\qquad if $\sigma\in {}^{\geqslant 2}D_{n-1,k}$
($\hat\sigma\in {}^{\geqslant 2}D_{n-1,k}$), then
\begin{equation}\label{i7}
\sigma_1 n\sigma_2\cdots\sigma_{n-1}\in {}^{\geqslant 2}D_{n,k+1},\qquad
\bar\sigma_1\bar n\sigma_2\cdots\sigma_{n-1}\in {}^{0,1}D_{n,k+2};
\end{equation}
\qquad if $\sigma\in {}^{0,\geqslant 2}D_{n-1,k}$
($\hat\sigma\in {}^1D_{n-1,k}$), then
\begin{equation}\label{i8}
\sigma_1 n\sigma_2\cdots\sigma_{n-1}\in {}^{\geqslant 2}D_{n,k},\qquad
\bar\sigma_1\bar n\sigma_2\cdots\sigma_{n-1}\in {}^{0,1}D_{n,k+1}.
\end{equation}
At position $i$ ($2\leqslant i\leqslant n-2$):

\qquad if $\sigma\in {}^*D_{n-1,k}$, then
\begin{equation}\label{i9}
\sigma_1\sigma_2\cdots\sigma_i n\sigma_{i+1}\cdots\sigma_{n-1}\in
\left\{
\begin{matrix}
{}^*D_{n,k} &\textrm{if $i\in\des(\sigma)$,}\\
{}^*D_{n,k+1} &\textrm{if $i\not\in\des(\sigma)$;}
\end{matrix}
\right.
\end{equation}
\qquad if $\hat\sigma\in {}^*D_{n-1,k}$, then
\begin{equation}\label{i10}
\bar\sigma_1\sigma_2\cdots\sigma_i\bar n\sigma_{i+1}\cdots\sigma_{n-1}\in
\left\{
\begin{matrix}
{}^*D_{n,k} &\textrm{if $i\in\des(\sigma)$,}\\
{}^*D_{n,k+1} &\textrm{if $i\not\in\des(\sigma)$;}
\end{matrix}
\right.
\end{equation}
At position $n-1$:
\begin{equation}\label{i11}
\sigma_1\sigma_2\cdots\sigma_{n-1}n\in {}^*D_{n,k};
\end{equation}
\begin{equation}\label{i12}
\bar\sigma_1\sigma_2\cdots\sigma_{n-1}\bar n\in {}^*D_{n,k+1}.
\end{equation}

\bigskip
Define now the {\it sub-Eulerian numbers} $D_{n,k}^1$, $D_{n,k}^{0,1}$,
$D_{n,k}^{\geqslant 2}$, $D_{n,k}^{0,\geqslant 2}$ by
$$
\begin{array}{rl}
D_{n,k}^1=&\#{}^1D_{n,k},\\
\noalign{\medskip}
D_{n,k}^{0,1}=&\#{}^{0,1}D_{n,k},\\
\noalign{\medskip}
D_{n,k}^{\geqslant 2}=&\#{}^{\geqslant 2}D_{n,k},\\
\noalign{\medskip}
D_{n,k}^{0,\geqslant 2}=&\#{}^{0,\geqslant 2}D_{n,k}.
\end{array}
$$
The sub-Eulerian numbers defined are obviously refinements of the
Eulerian numbers $D_{n,k}$,
and they offer a more accurate description of the descent distribution
of elements of $D_n$.

\bigskip
Since the map $\sigma\to\hat\sigma$ is an injection,
Lemma~\ref{lem31} then says that 
$$
\begin{array}{rl}
\#{}^1\bar D_{n,k}=&D_{n,k}^{0,\geqslant 2},\\
\noalign{\medskip}
\#{}^{0,1}\bar D_{n,k}=&D_{n,k}^{0,1},\\
\noalign{\medskip}
\#{}^{\geqslant 2}\bar D_{n,k}=&D_{n,k}^{\geqslant 2},\\
\noalign{\medskip}
\#{}^{0,\geqslant 2}\bar D_{n,k}=&D_{n,k}^1.
\end{array}
$$
Note that, in the insertion process above, the descent numbers can go up by
$2$ or go down by $1$, which never occurs in the case of $\mathfrak{S}_n$ or
$B_n$.
However, the introduction of sub-Eulerian numbers helps capture these changes.

\begin{Prop}\label{rec}
For $n\geqslant 3$,
the sub-Eulerian numbers satisfy the following recurrence relations
\begin{enumerate}
\item[\rm (i)]
$D_{n,k}^1=(2k-1)D_{n-1,k}^1+2(n-k)D_{n-1,k-1}^1+D_{n-1,k}^{0,1}
+D_{n-1,k-1}^{\geqslant 2}+D_{n-1,k}^{0,\geqslant 2}$,
\item[\rm (ii)]
$D_{n,k}^{0,1}=2(k-1)D_{n-1,k}^{0,1}+(2n-2k+1)D_{n-1,k-1}^{0,1}
+D_{n-1,k-1}^1+D_{n-1,k-2}^{\geqslant 2}+D_{n-1,k-1}^{0,\geqslant 2}$,
\item[\rm (iii)]
$D_{n,k}^{\geqslant 2}=(2k+1)D_{n-1,k}^{\geqslant 2}
+2(n-k-1)D_{n-1,k-1}^{\geqslant 2}+D_{n-1,k}^1+D_{n-1,k+1}^{0,1}
+D_{n-1,k}^{0,\geqslant 2}$,
\item[\rm (iv)]
$D_{n,k}^{0,\geqslant 2}=2kD_{n-1,k}^{0,\geqslant 2}
+(2n-2k-1)D_{n-1,k-1}^{0,\geqslant 2}+D_{n-1,k-1}^1+D_{n-1,k}^{0,1}
+D_{n-1,k-1}^{\geqslant 2}$.
\end{enumerate}
\end{Prop}
\begin{proof}
We only prove (i); the remaining assertions follow from similar considerations.
Elements of ${}^1D_{n,k}$ can be obtained by inserting $n$
$$
\begin{array}{rllcll}
\textrm{to}& \sigma\in {}^1D_{n-1,k-1}& \textrm{as in}&
(\ref{i1})& \textrm{with multiplicity}& D_{n-1,k-1}^1,\\
& \sigma\in {}^{0,1}D_{n-1,k}& & (\ref{i2})& & D_{n-1,k}^{0,1}\\
& \sigma\in {}^{\geqslant 2}D_{n-1,k-1}& & (\ref{i3})& &
D_{n-1,k-1}^{\geqslant 2}\\
& \sigma\in {}^{0,\geqslant 2}D_{n-1,k}& & (\ref{i4})& &
D_{n-1,k}^{0,\geqslant 2}\\
& \sigma\in {}^1D_{n-1,k}& & (\ref{i9})& & (k-1)D_{n-1,k}^1\\
& \sigma\in {}^1D_{n-1,k}& & (\ref{i9})& & [(n-2)-(k-1)]D_{n-1,k-1}^1\\
& \sigma\in {}^1D_{n-1,k}& & (\ref{i11})& & D_{n-1,k}^1
\end{array}
$$
and by inserting $-n$
$$
\begin{array}{rllcll}
\textrm{to}& \hat\sigma\in {}^1\bar D_{n-1,k}& \textrm{as in}&
(\ref{i10})& \textrm{with multiplicity}& (k-1)D_{n-1,k}^1,\\
& \hat\sigma\in {}^1\bar D_{n-1,k-1}& & (\ref{i10})& & 
[(n-2)-(k-1)]D_{n-1,k-1}^1,\\
& \hat\sigma\in {}^1\bar D_{n-1,k-1}& & (\ref{i12})& &D_{n-1,k-1}^1.
\end{array}
$$
Summing the multiplicities, we finally have
$$
D_{n,k}^1=(2k-1)D_{n-1,k}^1+2(n-k)D_{n-1,k-1}^1+D_{n-1,k}^{0,1}
+D_{n-1,k-1}^{\geqslant 2}+D_{n-1,k}^{0,\geqslant 2},
$$
as desired.
\end{proof}

It is convenient to record the sub-Eulerian numbers by generating functions.
For $n\geqslant 2$,
define the {\it sub-Eulerian polynomials} $D_n^1(t)$, $D_n^{0,1}(t)$,
$D_n^{\geqslant 2}(t)$, and $D_n^{0,\geqslant 2}(t)$ by
$$
\begin{array}{rl}
D_n^1(t)=&\displaystyle\sum_{k\geqslant 0}D_{n,k}^1t^k,\\
\noalign{\medskip}
D_n^{0,1}(t)=&\displaystyle\sum_{k\geqslant 0}D_{n,k}^{0,1}t^k,\\
\noalign{\medskip}
D_n^{\geqslant 2}(t)=&\displaystyle\sum_{k\geqslant 0}
D_{n,k}^{\geqslant 2}t^k,\\
\noalign{\medskip}
D_n^{0,\geqslant 2}(t)=&\displaystyle\sum_{k\geqslant 0}
D_{n,k}^{0,\geqslant 2}t^k.
\end{array}
$$
The first few values of $D_n^1(t)$, etc.,
are given in Tables~\ref{tab1}--\ref{tab4}.

\bigskip
It is clear that the Eulerian polynomial $D_n(t)$, $n\geqslant 2$,
satisfies
\begin{equation}\label{par}
D_n(t)=D_n^1(t)+D_n^{0,1}(t)+D_n^{\geqslant 2}(t)+D_n^{0,\geqslant 2}(t).
\end{equation}
The Eulerian polynomials of type $A$ and $B$ satisfy certain
difference-differential equations.
The same is true of the sub-Eulerian polynomials.

\begin{Prop}\label{dd}
For $n\geqslant 3$,
the sub-Eulerian polynomials satisfy the following difference-differential
equations:
\begin{enumerate}
\item[\rm (i)]
$D_n^1(t)=[2(n-1)t-1]D_{n-1}^1(t)+2t(1-t)(D_{n-1}^1)'(t)+D_{n-1}^{0,1}(t)
+tD_{n-1}^{\geqslant 2}(t)+D_{n-1}^{0,\geqslant 2}(t)$,
\item[\rm (ii)]
$D_n^{0,1}(t)=[(2n-1)t-2]D_{n-1}^{0,1}(t)+2t(1-t)(D_{n-1}^{0,1})'(t)
+tD_{n-1}^1(t)+t^2D_{n-1}^{\geqslant 2}(t)+tD_{n-1}^{0,\geqslant 2}(t)$,
\item[\rm (iii)]
$D_n^{\geqslant 2}(t)=[2(n-2)t+1]D_{n-1}^{\geqslant 2}(t)
+2t(1-t)(D_{n-1}^{\geqslant 2})'(t)+D_{n-1}^1(t)+t^{-1}D_{n-1}^{0,1}(t)
+D_{n-1}^{0,\geqslant 2}(t)$,
\item[\rm (iv)]
$D_n^{0,\geqslant 2}(t)=(2n-3)tD_{n-1}^{0,\geqslant 2}(t)
+2t(1-t)(D_{n-1}^{0,\geqslant 2})'(t)+tD_{n-1}^1(t)+D_{n-1}^{0,1}(t)
+tD_{n-1}^{\geqslant 2}(t)$.
\end{enumerate}
\end{Prop}
\begin{proof}
We only prove (i); the remaining assertions follow from similar reasoning.
Multiplying Proposition~\ref{rec} (i) by $t^k$ and summing over $k$ yields
\begin{equation}\label{dd1}
\begin{array}{rl}
\displaystyle\sum_{k=1}^n D_{n,k}^1t^k
=&\displaystyle\sum_{k=1}^n(2k-1)D_{n-1,k}^1t^k
+\sum_{k=1}^n 2(n-k)D_{n-1,k-1}^1t^k+\sum_{k=1}^n D_{n-1,k}^{0,1}t^k\\
&\displaystyle
\qquad+\sum_{k=1}^n D_{n-1,k-1}^{\geqslant 2}t^k
+\sum_{k=1}^n D_{n-1,k}^{0,\geqslant 2}t^k\\
\noalign{\medskip}
=&I+II+III+IV+V.
\end{array}
\end{equation}
The left hand side of (\ref{dd1}) is equal to $D_n^1(t)$ because
$D_{n,0}^1=0$, while terms on the right hand side are equal respectively
to
$$
\begin{array}{rl}
I=&\displaystyle
2t\sum_{k=1}^n kD_{n-1,k}^1t^{k-1}-\sum_{k=1}^n D_{n-1,k}^1t^k
=2t(D_{n-1}^1)'(t)-D_{n-1}^1(t),\\
\noalign{\medskip}
II=&\displaystyle\sum_{k=0}^{n-1}2(n-k-1)D_{n-1,k}^1t^{k+1}
=2(n-1)t\sum_{k=0}^{n-1}D_{n-1,k}^1t^k
-2t^2\sum_{k=0}^{n-1}kD_{n-1,k}^1t^{k-1}\\
=&2(n-1)tD_{n-1}^1(t)-2t^2(D_{n-1}^1)'(t),\\
\noalign{\medskip}
III=&\displaystyle\sum_{k=0}^{n-1}D_{n-1,k}^{0,1}t^k=D_{n-1}^{0,1}(t),\\
\noalign{\medskip}
IV=&\displaystyle\sum_{k=0}^{n-1}D_{n-1,k}^{\geqslant 2}t^{k+1}
=tD_{n-1}^{\geqslant 2}(t),\\
\noalign{\medskip}
V=&D_{n-1}^{0,\geqslant 2}(t),
\end{array}
$$
because $D_{n,0}^1=D_{n-1,n}^1=D_{n,0}^{0,1}=D_{n-1,n}^{0,1}=
D_{n-1,0}^{0,\geqslant 2}=D_{n-1,n}^{0,\geqslant 2}=0$.
Hence, (i) follows.
\end{proof}

\bigskip
In view of the decomposition property (\ref{par}) of the sub-Eulerian
polynomials, the following is the immediate consequence of
Proposition~\ref{dd}.

\begin{Cor}\label{dnrec}
We have
$$
\begin{array}{rl}
D_n(t)=&2ntD_{n-1}(t)+2t(1-t)D_{n-1}'(t)+(t^{-1}-t)D_{n-1}^{0,1}(t)
+(1-t)^2D_{n-1}^{\geqslant 2}(t)\\
\noalign{\medskip}
&\qquad+2(1-t)D_{n-1}^{0,\geqslant 2}(t).
\end{array}
$$
\end{Cor}
\begin{proof}
Summing Proposition~\ref{dd} (i)--(iv) and using (\ref{par}).
\end{proof}

The above recurrence relation will be revisited later when further properties
of the sub-Eulerian polynomials are found.

\newpage
\begin{table}
\begin{center}
\begin{tabular}[t]{c @{\qquad} l}
$n$& $D_n^1(t)$\\
\hline
$2$& $t$\\
$3$& $3t+3t^2$\\
$4$& $7t+34t^2+7t^3$\\
$5$& $15t+225t^2+225t^3+15t^4$\\
$6$& $31t+1196t^2+3306t^3+1196t^4+31t^5$\\
\hline
\end{tabular}
\end{center}
\caption{The sub-Eulerian polynomial $D_n^1(t)$ for $n=2,\ldots,6$.}
\label{tab1}
\end{table}

\bigskip
\begin{table}
\begin{center}
\begin{tabular}[t]{c @{\qquad} l}
$n$& $D_n^{0,1}(t)$\\
\hline
$2$& $t^2$\\
$3$& $5t^2+t^3$\\
$4$& $17t^2+30t^3+t^4$\\
$5$& $49t^2+303t^3+127t^4+t^5$\\
$6$& $129t^2+2132t^3+3030t^4+468t^5+t^6$\\
\hline
\end{tabular}
\end{center}
\caption{The sub-Eulerian polynomial $D_n^{0,1}(t)$ for $n=2,\ldots,6$.}
\label{tab2}
\end{table}

\bigskip
\begin{table}
\begin{center}
\begin{tabular}[t]{c @{\qquad} l}
$n$& $D_n^{\geqslant 2}(t)$\\
\hline
$2$& $1$\\
$3$& $1+5t$\\
$4$& $1+30t+17t^2$\\
$5$& $1+127t+303t^2+49t^3$\\
$6$& $1+468t+3030t^2+2132t^3+129t^4$\\
\hline
\end{tabular}
\end{center}
\caption{The sub-Eulerian polynomial $D_n^{\geqslant 2}(t)$
for $n=2,\ldots,6$.}\label{tab3}
\end{table}

\bigskip
\begin{table}
\begin{center}
\begin{tabular}[t]{c @{\qquad} l}
$n$& $D_n^{0,\geqslant 2}(t)$\\
\hline
$2$& $t$\\
$3$& $3t+3t^2$\\
$4$& $7t+34t^2+7t^3$\\
$5$& $15t+225t^2+225t^3+15t^4$\\
$6$& $31t+1196t^2+3306t^3+1196t^4+31t^5$\\
\hline
\end{tabular}
\end{center}
\caption{The sub-Eulerian polynomial $D_n^{0,\geqslant 2}(t)$
for $n=2,\ldots,6$.}\label{tab4}
\end{table}

\newpage
\section{Generating Functions}

Proposition~\ref{dd} offers an efficient way of computing $D_n^1(t)$, etc.
This section discusses a convenient way of recording them, e.g.,
by their exponential generating functions.
Define the exponential generating functions for the sub-Eulerian polynomials
by
\begin{equation}\label{egf}
\begin{array}{c}
D^1(x,t)=\displaystyle\sum_{n\geqslant 2}D_n^1(t)\frac{x^n}{n!},\\
D^{0,1}(x,t)=\displaystyle\sum_{n\geqslant 2}D_n^{0,1}(t)\frac{x^n}{n!},\\
D^{\geqslant 2}(x,t)=\displaystyle
\sum_{n\geqslant 2}D_n^{\geqslant 2}(t)\frac{x^n}{n!},\\
D^{0,\geqslant 2}(x,t)=\displaystyle
\sum_{n\geqslant 2}D_n^{0,\geqslant 2}(t)\frac{x^n}{n!}.
\end{array}
\end{equation}
Here, we are interested in obtaining closed form expressions for the
right hand sides of (\ref{egf}).

\bigskip
Denote by $\mathfrak{S}_n(t)$ and $B_n(t)$ the $n$-th Eulerian polynomials of
type $A$ and $B$ for $\mathfrak{S}_n$ and $B_n$, respectively.
In the case of $\mathfrak{S}(x,t)=\sum_{n\geqslant 0}\mathfrak{S}_n(t)x^n/n!$,
one determines $f_{n,k}$ for which $\mathfrak{S}_n(t)/(1-t)^{n+1}=
\sum_{k\geqslant 0}f_{n,k}t^k$, then multiply by $x^n/n!$, sum over $n$,
and replace $x$ by $x(1-t)$.
This same procedure is complicated for $D^1(x,t)$, etc.,
because they are coupled via Proposition~\ref{dd} (i)--(iv).
Also, it is not clear at present what are the right denominators $Q(t)$
for the rational generating function $D_n^*(t)/Q(t)$.

\bigskip
An alternative, which we believe to be new,
way of computing them is as follows.
Let us show how it works in the case of $\mathfrak{S}(x,t)$.
By taking partial derivatives of $\mathfrak{S}(x,t)$ with respect to $x$
and $t$,
and making use of the difference-differential equation which $\mathfrak{S}_n(t)$
satisfy, namely,
\begin{equation}
\mathfrak{S}_n(t)=[(n-1)t+1]\mathfrak{S}_{n-1}(t)+t(1-t)\mathfrak{S}_{n-1}'(t),
\end{equation}
we obtain that $\mathfrak{S}=\mathfrak{S}(x,t)$ satisfies the following first
order linear partial differential equation
\begin{equation}\label{pdea}
t(t-1)\mathfrak{S}_t+(1-xt)\mathfrak{S}_x=\mathfrak{S},
\end{equation}
which, together with the initial condition $\mathfrak{S}(0,t)=1$,
uniquely determine $\mathfrak{S}$, that is,
\begin{equation}
\mathfrak{S}(x,t)=\frac{(1-t)e^{x(1-t)}}{1-te^{x(1-t)}}.
\end{equation}
For completeness, we record the PDE which $B(x,t)=\sum_{n\geqslant 0}B_n(t)
x^n/n!$ satisfies, namely
\begin{equation}\label{b}
2t(t-1)B_t+(1-2xt)B_x=(1+t)B.
\end{equation}
The concerned initial condition is $B(0,t)=1$,
and the solution to (\ref{b}) is of course
\begin{equation}
B(x,t)=\frac{(1-t)e^{x(1-t)}}{1-te^{2x(1-t)}}.
\end{equation}
We shall not give the details of the proof of (\ref{pdea}) and
(\ref{b}) but encourage interested readers to work them out by imitating
the proof of Proposition~\ref{pde} below.

\begin{Prop}\label{pde}
We have
\begin{enumerate}
\item[\rm (i)]
$2t(t-1)D_t^1+(1-2xt)D_x^1=-D^1+D^{0,1}+tD^{\geqslant 2}+D^{0,\geqslant 2}
+xt$;
\item[\rm (ii)]
$2t(t-1)D_t^{0,1}+(1-2xt)D_x^{0,1}=tD^1+(t-2)D^{0,1}+t^2D^{\geqslant 2}
+tD^{0,\geqslant 2}+xt^2$;
\item[\rm (iii)]
$2t(t-1)D_t^{\geqslant 2}+(1-2xt)D_x^{\geqslant 2}=D^1+t^{-1}D^{0,1}
+(1-2t)D^{\geqslant 2}+D^{0,\geqslant 2}+x$;
\item[\rm (iv)]
$2t(t-1)D_t^{0,\geqslant 2}+(1-2xt)D_x^{0,\geqslant 2}=tD^1+D^{0,1}
+tD^{\geqslant 2}-tD^{0,\geqslant 2}+xt$.
\end{enumerate}
\end{Prop}
\begin{proof}
We only prove (iii); the remaining assertions follow from similar reasoning.
We have
$$
\begin{array}{rl}
2t(1-t)D_t^{\geqslant 2}=&\displaystyle
\sum_{n\geqslant 2}2t(1-t)(D_n^{\geqslant 2})'(t)\frac{x^n}{n!}\\
\noalign{\medskip}
=&\displaystyle
\sum_{n\geqslant 2}\{D_{n+1}^{\geqslant 2}(t)-[2nt-(2t-1)]D_n^{\geqslant 2}(t)
-D_n^{0,\geqslant 2}(t)-D_n^1(t)\\
&\displaystyle\quad-t^{-1}D_n^{0,1}(t)\}\frac{x^n}{n!}\\
\noalign{\medskip}
=&\displaystyle
\sum_{n\geqslant 2}D_{n+1}^{\geqslant 2}(t)\frac{x^n}{n!}
-2xt\sum_{n\geqslant 2}D_n^{\geqslant 2}(t)\frac{x^{n-1}}{(n-1)!}\\
&\displaystyle\quad
+(2t-1)\sum_{n\geqslant 2}D_n^{n\geqslant 2}(t)\frac{x^n}{n!}
-\sum_{n\geqslant 2}D_n^{0,\geqslant 2}(t)\frac{x^n}{n!}\\
&\displaystyle\quad
-\sum_{n\geqslant 2}D_n^1(t)\frac{x^n}{n!}
-t^{-1}\sum_{n\geqslant 2}D_n^{0,1}(t)\frac{x^n}{n!}\\
\noalign{\medskip}
=&\displaystyle
D_x^{\geqslant 2}-D_2^{\geqslant 2}(t)x-2xtD_x^{\geqslant 2}
+(2t-1)D^{\geqslant 2}-D^{0,\geqslant 2}-D^1\\
&\quad-t^{-1}D^{0,1},
\end{array}
$$
from which (iii) follows, because $D_2^{\geqslant 2}(t)=1$.
\end{proof}

\bigskip
The above set of PDEs are readily solved by the method of characteristics
\cite{ch89}. The idea is to solve for characteristic $x=x(t;A)$,
which satisfies $dx/dt=(1-2xt)/2t(t-1)$ and is parametrized by $A$,
where $A$ is the constant of integration;
for $x$ and $t$ related by $x=x(t;A)$,
the partial derivatives become total derivatives with respect to $t$,
and the PDEs are thus reduced to ODEs and the solutions of which solve
the given PDEs. (Here, the constants of integration of the latter ODEs
are arbitrary functions of $A$, and whose forms can be determined by the
initial conditions.)

\begin{Th}\label{sepegf}
We have
\begin{enumerate}
\item[\rm (i)]
$D^{0,\geqslant 2}(x,t)=\displaystyle
\frac{t(e^{x(1-t)}-1)^2}{2(1-t)(1-te^{2x(1-t)})}=D^1(x,t)$,
\item[\rm (ii)]
$D^{\geqslant 2}(x,t)=\displaystyle
\frac{-1-x(1-t)+e^{x(1-t)}}{(1-t)(1-te^{2x(1-t)})}$,
\item[\rm (iii)]
$D^{0,1}(x,t)=\displaystyle
\frac{t^2e^{x(1-t)}-t^2[1-x(1-t)]e^{2x(1-t)}}{(1-t)(1-te^{2x(1-t)})}$.
\end{enumerate}
\end{Th}
\begin{proof}
Let $x=x(t)$ be such that $dx/dt=(1-2xt)/2t(t-1)$.
The latter linear first order ODE is readily solved, e.g.,
$t^{1/2}e^{x(1-t)}=A$, where $A$ is a constant.
On $x=x(t)$, the left hand sides of Proposition~\ref{pde} become total
derivatives with respect to $t$.
In particular, Proposition~\ref{pde} (iv) is then
\begin{equation}\label{ode1}
2t(t-1)\frac{dD^{0,\geqslant 2}}{dt}=tD^1+D^{0,1}
+tD^{\geqslant 2}-tD^{0,\geqslant 2}.
\end{equation}
Differentiating (\ref{ode1}) with respect to $t$, using
$dx/dt=(1-2xt)/2t(t-1)$, multiplying by $2t(t-1)$,
followed by substituting the $t$-derivatives
on the right by the corresponding right hand side in Proposition~\ref{pde},
we have that $D^{0,\geqslant 2}$ satisfies the following ODE,
\begin{equation}\label{ode2}
4t(t-1)^2\frac{d^2D^{0,\geqslant 2}}{dt^2}+2(t-1)(3t-1)\frac{dD^{0,\geqslant
2}}{dt}-4D^{0,\geqslant 2}=1,
\end{equation}
whose solution (obtained by {\it Mathematica}) is
\begin{equation}\label{ode3}
D^{0,\geqslant 2}(x,t)=\frac{1-t+4(1+t)C_1(A)-8it^{1/2}C_2(A)}{4(t-1)}.
\end{equation}
Here, $C_1$, and $C_2$ are arbitrary functions of $A$.
The initial condition $D^{0,\geqslant 2}(0,t)=0$ implies that
$$
\frac{1-t+4(1+t)C_1(t^{1/2})-8it^{1/2}C_2(t^{1/2})}{4(t-1)}=0
$$
so that
$$
C_2(t)=\frac{1-t^2+4(1+t^2)C_1(t)}{8it}.
$$
Reporting this back in (\ref{ode3}) and replacing $A$ by $t^{1/2}e^{x(1-t)}$
yields that
\begin{equation}\label{d02}
D^{0,\geqslant 2}(x,t)=\frac{(e^{x(1-t)}-1)[te^{x(1-t)}+1+4(1-te^{x(1-t)})
C_1(t^{1/2}e^{x(1-t)})]}{4(t-1)e^{x(1-t)}}.
\end{equation}
The initial condition $D_x^{0,\geqslant 2}(0,t)=0$ then yields that
$$
0=\frac{(1-t)(t+1+4(1-t)C_1(t^{1/2}))}{4(t-1)},
$$
from which we deduce that
$$
C_1(t)=\frac{t^2+1}{4(t^2-1)}.
$$
Substituting $C_1(t)$ back in (\ref{d02}) and after some algebra,
the first equality in (i) follows.
To prove the second equality in (i),
we differentiate
\begin{equation}
2t(t-1)\frac{dD^1}{dt}=-D^1+D^{0,1}+tD^{\geqslant 2}+D^{0,\geqslant 2}
+xt
\end{equation}
with respect to $t$ along the characteristic $t^{1/2}e^{x(1-t)}=A$,
multiply the resulting equation by $2t(t-1)$,
and then substitute the $t$-derivatives by the corresponding right hand
sides in Proposition~\ref{pde}, the end result being
\begin{equation}\label{ode4}
4t(t-1)^2\frac{d^2D^1}{dt^2}+2(t-1)(3t-1)\frac{dD^1}{dt}-4D^1=1.
\end{equation}
Note that (\ref{ode4}) is the same as (\ref{ode2}).
Since the initial conditions $D^1(0,t)=D_x^1(0,t)=0$ are also the same as
those for $D^{0,\geqslant 2}$,
we conclude that the second equality in (i) holds.

We have, from Proposition~\ref{pde} (iii) and (iv), that
$$
t(t-1)\frac{dD^{\geqslant 2}}{dt}+tD^{\geqslant 2}
=(t-1)\frac{dD^{0,\geqslant 2}}{dt}+D^{0,\geqslant 2},
$$
which can be written simply as
\begin{equation}\label{ode5}
t\frac{d}{dt}((t-1)D^{\geqslant 2})=\frac{d}{dt}((t-1)D^{0,\geqslant 2}).
\end{equation}
The right hand side of (\ref{ode5}) is readily computed, i.e.,
$$
\begin{array}{rl}
\displaystyle\frac{d}{dt}((t-1)D^{0,\geqslant 2})=&\displaystyle
\frac{d}{dt}\left(\frac{t(e^{x(1-t)}-1)^2}{2(te^{2x(1-t)}-1)}\right)
=\frac{d}{dt}\left(\frac{A^2-2t^{1/2}A+t}{2(A^2-1)}\right)\\
\noalign{\medskip}
=&\displaystyle\frac{-t^{-1/2}A+1}{2(A^2-1)},
\end{array}
$$
so that (\ref{ode5}) now reads
$$
t\frac{d}{dt}((t-1)D^{\geqslant 2})=\frac{-t^{-1/2}A+1}{2(A^2-1)}.
$$
Straightforward integration yields that
$$
\begin{array}{rl}
D^{\geqslant 2}=&\displaystyle
\frac{2t^{-1/2}A+\ln t}{2(1-t)(1-A^2)}+\frac{f(A)}{t-1}\\
\noalign{\medskip}
=&\displaystyle
\frac{2e^{x(1-t)}+\ln t}{2(1-t)(1-te^{2x(1-t)})}+
\frac{f(t^{1/2}e^{x(1-t)})}{t-1}.
\end{array}
$$
The condition $D^{\geqslant 2}(0,t)=0$ then forces
$$
f(t)=\frac{1+\ln t}{1-t^2}.
$$
With this $f$, and after some algebra, (ii) follows.

We have, from Proposition~\ref{pde} (i) and (ii), that
\begin{equation}\label{ode6}
t(t-1)\frac{dD^{0,1}}{dt}+D^{0,1}=t^2(t-1)\frac{dD^1}{dt}+tD^1,
\end{equation}
By (i),
$$
D^1=\frac{t(e^{x(1-t)}-1)^2}{2(t-1)(te^{2x(1-t)}-1)}=
\frac{A^2-2t^{1/2}A+t}{2(t-1)(A^2-1)}.
$$
The right hand side of (\ref{ode6}) is easily computed, i.e.,
$$
t^2(t-1)\frac{dD^1}{dt}+tD^1=\frac{At^{3/2}-A^2t}{2(t-1)(A^2-1)},
$$
so that (\ref{ode6}) can be rewritten as
$$
\frac{dD^{0,1}}{dt}+\frac{D^{0,1}}{t(t-1)}=
\frac{At^{1/2}-A^2}{2(t-1)(A^2-1)}.
$$
By the standard procedure for solving first order linear ODEs,
we have that
\begin{equation}\label{ode7}
\begin{array}{rl}
D^{0,1}=&\displaystyle\frac{t}{t-1}\left\{
\frac{2At^{1/2}-A^2\ln t}{2(A^2-1)}+g(A)\right\}\\
\noalign{\medskip}
=&\displaystyle
\frac{2t^2e^{x(1-t)}-t^2e^{2x(1-t)}\ln t}{2(1-t)(1-te^{2x(1-t)})}
+\frac{tg(t^{1/2}e^{x(1-t)})}{t-1},
\end{array}
\end{equation}
where $g$ is an arbitrary function.
The condition $D^{0,1}(0,t)=0$ then implies that
$$
g(t)=\frac{t^2\ln t-t^2}{t^2-1}.
$$
Reporting $g$ back in (\ref{ode7}), and after some algebra, (iii) follows.
\end{proof}

\bigskip
Some properties of $D_n^1(t)$, etc., are readily deduced from
Theorem~\ref{sepegf}.

\begin{Cor}\label{refl}
For $n\geqslant 2$,
\begin{enumerate}
\item[\rm (i)]
$D_n^1(t)=D_n^{0,\geqslant 2}(t)$,
\item[\rm (ii)]
$D_n^{0,1}(t)=t^nD_n^{\geqslant 2}(t^{-1})$,
\item[\rm (iii)]
$D_n^{\geqslant 2}(t)=t^nD_n^{0,1}(t^{-1})$,
\item[\rm (iv)]
$D_n^{0,\geqslant 2}(t)=t^nD_n^{0,\geqslant 2}(t^{-1})$,
\item[\rm (v)]
$D_n^1(t)=t^nD_n^1(t^{-1})$,
\end{enumerate}
\end{Cor}
\begin{proof}
Since $D^1(x,t)=D^{0,\geqslant 2}(x,t)$,
equating the coefficients of $x^n$ of $D^1(x,t)$ and
$D^{0,\geqslant 2}(x,t)$, (i) follows.
Since
$$
\begin{array}{rl}
D^{\geqslant 2}(xt,t^{-1})=&\displaystyle\
\frac{-1-xt(1-t^{-1})+e^{x(1-t^{-1})}}{(1-t^{-1})(1-t^{-1}e^{2xt(1-t^{-1})})}\\
\noalign{\medskip}
=&\displaystyle
\frac{t^2e^{2x(1-t)}[-1+x(1-t)+e^{-x(1-t)}]}{(1-t)(1-te^{2x(1-t)})}\\
\noalign{\medskip}
=&\displaystyle
\frac{t^2e^{x(1-t)}-t^2[1-x(1-t)]e^{2x(1-t)}}{(1-t)(1-te^{2x(1-t)})}
=D^{0,1}(x,t),
\end{array}
$$
equating the coefficients of $x^n$ on both sides, (ii) follows.
To prove (iii), we replace $t$ by $t^{-1}$, and $x$ by $xt$
in $D^{\geqslant 2}(xt,t^{-1})=D^{0,1}(x,t)$,
the end result being
$$
D^{\geqslant 2}(x,t)=D^{\geqslant 2}((xt)(t^{-1}),(t^{-1})^{-1})
=D^{0,1}(xt,t^{-1}),
$$
from which (iii) follows.
By a calculation similar to that in (ii),
we have that $D^{0,\geqslant 2}(xt,t^{-1})=D^{0,\geqslant 2}(x,t)$
from which (iv) follows.
(v) follows from (i) and (iv).
\end{proof}

\begin{Rem}{\rm
Corollary~\ref{refl} (ii)--(v) state that $D_n^1(t)$, etc., have certain
reflectional symmetry.
For example, the polynomial $t^nD_n^{\geqslant 2}(t^{-1})$ is obtained
by reflecting the coefficients of $D_n^{\geqslant 2}(t)$ about $x^{n/2}$.
Thus, Corollary~\ref{refl} (ii)--(iii) (resp., (iv)--(v)) say that
$D_n^{0,1}(t)$ (resp., $D_n^1(t)$) are obtained from $D_n^{\geqslant 2}(t)$
(resp., $D_n^{0,\geqslant 2}(t)$) by reflection, and vice versa.
Corollary~\ref{refl} (i) further says that $D_n^1(t)$ and
$D_n^{0,\geqslant 2}(t)$ are themselves reflectionally symmetric.
However, $D_n^{0,1}(t)$ and $D_n^{\geqslant 2}(t)$ do not enjoy this
property.
}\end{Rem}

Corollary~\ref{refl} (i) enables a furthur simplification of the recurrence
relation for $D_n(t)$ derived in the previous section.

\begin{Cor}
We have
\begin{equation}\label{dnrec2}
\begin{array}{rl}
D_n(t)=&[(2n-1)t+1]D_{n-1}(t)+2t(1-t)D_{n-1}'(t)\\
\noalign{\medskip}
&\qquad+(1-t)[t^{-1}D_{n-1}^{0,1}(t)-tD_{n-1}^{\geqslant 2}(t)].
\end{array}
\end{equation}
\end{Cor}
\begin{proof}
By Corollary~\ref{refl} (i) and (\ref{par}),
the last three terms on the right hand side of Corollary~\ref{dnrec}
can be written as
$$
\begin{array}{rl}
&(t^{-1}-t)D_{n-1}^{0,1}(t)+(1-t)^2D_{n-1}^{\geqslant 2}(t)
+2(1-t)D_{n-1}^{0,\geqslant 2}(t)\\
=&(1-t)[(1+t^{-1})D_{n-1}^{0,1}(t)+(1-t)D_{n-1}^{\geqslant 2}(t)
+D_{n-1}^1(t)+D_{n-1}^{0,\geqslant 2}(t)]\\
=&(1-t)D_{n-1}(t)+(1-t)[t^{-1}D_{n-1}^{0,1}(t)-tD_{n-1}^{\geqslant 2}(t)].
\end{array}
$$
Aftering regrouping terms, the corollary follows.
\end{proof}

\begin{Rem}{\rm
The first line of (\ref{dnrec2})
is precisely the recurrence relation for the Eulerian polynomial $B_n(t)$
of type $B$.
}\end{Rem}

\bigskip
Denote by $D(x,t)$ the exponential generating function for the Eulerian
polynomial $D_n(t)$ of type $D$, i.e.,
\begin{equation}
D(x,t)=\sum_{n\geqslant 0}D_n(t)\frac{x^n}{n!},
\end{equation}
where $D_0(t)=1$.
Since $D_1(t)=1$, it is clear that
\begin{equation}\label{dxt}
D(x,t)=1+x+D^1(x,t)+D^{0,1}(x,t)+D^{\geqslant 2}(x,t)+D^{0,\geqslant 2}(x,t).
\end{equation}
Computing $D(x,t)$ by (\ref{dxt}) and Theorem~\ref{sepegf},
we obtain the following result of Brenti \cite{b94}.

\begin{Cor}
We have
\begin{equation}
D(x,t)=\frac{(1-t)e^{x(1-t)}-xt(1-t)e^{2x(1-t)}}{1-te^{2x(1-t)}}.
\end{equation}
\end{Cor}

\bigskip
We have now the exponential generating functions for the sub-Eulerian
polynomials at our disposal.
By reversing the procedure described in the paragraph following the
definitions of $D^1(x,t)$, etc., we have the following results.

\begin{Cor}
For $n\geqslant 2$,
\begin{enumerate}
\item[\rm (i)]
$\displaystyle\frac{D_n^{0,\geqslant 2}(t)}{t(1-t)^{n-1}}
=\sum_{k\geqslant 0}\{2^{n-1}[(k+1)^n+k^n]-(2k+1)^n\}t^k
=\frac{D_n^1(t)}{t(1-t)^{n-1}}$,
\item[\rm (ii)]
$\displaystyle\frac{D_n^{\geqslant 2}(t)}{(1-t)^{n-1}}
=\sum_{k\geqslant 0}[(2k+1)^n-(n+1)(2k)^n]t^k$,
\item[\rm (iii)]
$\displaystyle\frac{D_n^{0,1}(t)}{t^2(1-t)^{n-1}}
=\sum_{k\geqslant 0}[(2k+1)^n-2^n(k+1)^n+n2^{n-1}(k+1)^{n-1}]t^k$.
\end{enumerate}
\end{Cor}

\begin{proof}
We only prove the first equality in (i);
the remaining assertions follow from similar reasoning.
$$
\begin{array}{rl}
\displaystyle\sum_{n\geqslant 2}
\frac{D_n^{0,\geqslant 2}(t)x^n}{t(1-t)^{n-1}n!}
=&\displaystyle\frac{1-t}{t}D^{0,\geqslant 2}(\frac{x}{1-t},t)\\
\noalign{\medskip}
=&\displaystyle\frac{(e^x-1)^2}{2(1-te^{2x})}\\
\noalign{\medskip}
=&\displaystyle\frac{(e^{2x}-2e^x+1)}{2}\sum_{k\geqslant 0}t^ke^{2kx}\\
\noalign{\medskip}
=&\displaystyle\sum_{k\geqslant 0}t^k
\frac{(e^{2(k+1)x}-2e^{(2k+1)x}+e^{2kx})}{2}\\
\noalign{\medskip}
=&\displaystyle\sum_{k\geqslant 0}t^k\sum_{n\geqslant 0}
[2^{n-1}(k+1)^n-(2k+1)^n+2^{n-1}k^n]\frac{x^n}{n!}.
\end{array}
$$
Equating the coefficient of $x^n$ on both sides, the result follows.
\end{proof}

\section{A recurrence relation}

In the previous section, we computed the PDEs which the generating functions
of the sub-Eulerian polynomials satisfy.
By solving these PDEs, we obtained the generating functions.
It is important to note that the PDEs encode the recurrence relations
for the sub-Eulerian polynomials.
This procedure for computing generating functions can be reversed,
i.e., starting from the generating function,
we determine a PDE of minimal order which the generating function satisfies;
the desired recurrence relation then follows upon equating coefficients.

\bigskip
We shall apply this reversed procedure to $D(x,t)$ to obtain a recurrence
relation for $D_n(t)$.
Proposition~\ref{pde} and (\ref{dxt}) suggest the partial differential
operator $2t(t-1)\partial_t+(1-2xt)\partial_x$ for $D(x,t)$.
The proofs of the following lemma are easy, albeit tedious,
exercises of calculus which we omit.

\begin{Lem}\label{cal}
We have
\begin{enumerate}
\item[\rm(i)]
$[2t(t-1)\partial_t+(1-2xt)\partial_x]D(x,t)=\displaystyle
\frac{(1-t)[(1+t)e^{x(1-t)}-te^{2x(1-t)}]}{1-te^{2x(1-t)}}$;
\item[\rm(ii)]
$[2t(t-1)\partial_t+(1-2xt)\partial_x]^2D(x,t)=\displaystyle
\frac{(1-t)[(1+3t^2)e^{x(1-t)}-2t^2e^{2x(1-t)}]}{1-te^{2x(1-t)}}$;
\item[\rm (iii)]
$[2t(t-1)\partial_t+(1-2xt)\partial_x]^2=4t^2(1-t)^2\partial_{tt}
-4t(1-t)(1-2xt)\partial_{xt}+(1-2xt)^2\partial_{xx}
-2t(1-2x)\partial_x+4t(1-t)(1-2t)\partial_t$.
\end{enumerate}
\end{Lem}

\begin{Prop}
The function $D=D(x,t)$ satisfies the following PDE:
$$
\begin{array}{rl}
&4t^2(1-t)^2[x(1+t)-1]D_{tt}-
4t(1-t)(1-2xt)[x(1+t)-1]D_{xt}+\\
&(1-2xt)^2[x(1+t)-1]D_{xx}+2t(1-t)^2[-2+(3+t)x]D_t+\\
&[4t-x(1+3t)^2+2x^2t(3+2t+3t^2)]D_x+(1-t)^2D=0.
\end{array}
$$
\end{Prop}
\begin{proof}
Use Lemma~\ref{cal}(i)--(ii) to form the expression
\begin{equation}\label{express}
\begin{array}{l}
\alpha[2t(t-1)\partial_t+(1-2xt)\partial_x]^2D(x,t)
+\beta[2t(t-1)\partial_t+(1-2xt)\partial_x]D(x,t)\\
\qquad\qquad=\displaystyle\frac{(1-t)\{[\alpha(1+t)+\beta(1+3t^2)]e^{x(1-t)}
+(-\alpha t-2\beta t^2)e^{2x(1-t)}\}}{1-te^{2x(1-t)}}\\
\end{array}
\end{equation}
where $\alpha=\alpha(x,t)$, and $\beta=\beta(x,t)$ are functions to
be determined. Now assume that the right side of (\ref{express})
is equal to
\begin{equation}\label{rhs}
\frac{(1-t)[\gamma e^{x(1-t)}-\gamma xte^{2x(1-t)}]}{1-te^{2x(1-t)}}=
\gamma D(x,t),
\end{equation}
where $\gamma=\gamma(x,t)$ is a function to be determined.
Equating the coefficients of $e^{x(1-t)}$ and $e^{2x(1-t)}$
in the numerators of (\ref{express}) and (\ref{rhs}), we have that
$$
\alpha(1+t)+\beta(1+3t^2)=\gamma,\quad\textrm{and}\quad
-\alpha t-2\beta t^2=-\gamma xt,
$$
whose formal solutions are $\alpha=\gamma(1-t)^{-2}[x(1+3t^3)-2t]$,
and $\beta=-\gamma(1-t)^{-2}[x(1+t)-1]$.
Setting $\gamma=(1-t)^2$, we get $\alpha=[x(1+3t^2)-2t]$ and
$\beta=-[x(1+t)-1]$. Reporting these choices of $\alpha$, $\beta$, and
$\gamma$ in (\ref{express}), and using Lemma~\ref{cal}(iii),
the result follows.
\end{proof}

\bigskip
It is trivial to see that
$$
\begin{array}{l}
\displaystyle
D=\sum_{n\geqslant 0}D_n(t)\frac{x^n}{n!},\quad
D_t=\sum_{n\geqslant 0}D_n'(t)\frac{x^n}{n!},\quad
D_x=\sum_{n\geqslant 0}D_{n+1}(t)\frac{x^n}{n!},\\
\displaystyle
D_{tt}=\sum_{n\geqslant 0}D_n''(t)\frac{x^n}{n!},\quad
D_{xt}=\sum_{n\geqslant 0}D_{n+1}'(t)\frac{x^n}{n!},\quad
D_{xx}=\sum_{n\geqslant 0}D_{n+2}(t)\frac{x^n}{n!}.
\end{array}
$$
We are now ready to state the main result of this section.

\begin{Th}\label{dnrec5}
For $n\geqslant 1$,
the Eulerian polynomials $D_n(t)$ of type $D$ satisfy
\begin{equation}\label{dnrec3}
\begin{array}{rl}
D_{n+2}(t)=&[n(1+5t)+4t]D_{n+1}(t)\\
&+4t(1-t)D_{n+1}'(t)\\
&+[(1-t)^2-n(1+3t)^2-4n(n-1)t(1+2t)]D_n(t)\\
&-[4nt(1-t)(1+3t)+4t(1-t)^2]D_n'(t)\\
&-4t^2(1-t)^2D_n''(t)\\
&+[2n(n-1)t(3+2t+3t^2)-4n(n-1)(n-2)t^2(1+t)]D_{n-1}(t)\\
&+[2nt(1-t)^2(3+t)+8n(n-1)t^2(1-t)(1+t)]D_{n-1}'(t)\\
&+4nt^2(1-t)^2(1+t)D_{n-1}''(t).
\end{array}
\end{equation}
\end{Th}
\begin{proof}
Denote by $I,II,\ldots,VI$ the successive terms of the PDE as in
the previous proposition. Then
$$
\begin{array}{rl}
I=&\displaystyle
[-4t^2(1-t)^2+4t^2(1-t)^2(1+t)x]\sum_{n\geqslant 0}D_n''(t)\frac{x^n}{n!},\\
II=&\displaystyle
[4t(1-t)-4t(1-t)(1+3t)x+8t^2(1-t)(1+t)x^2]
\sum_{n\geqslant 0}D_{n+1}'(t)\frac{x^n}{n!},\\
III=&\displaystyle
[-1+(1+5t)x-4t(1+2t)x^2+4t^2(1+t)x^3]
\sum_{n\geqslant 0}D_{n+2}(t)\frac{x^n}{n!},\\
IV=&\displaystyle
[-4t(1-t)^2+2t(1-t)^2(3+t)x]
\sum_{n\geqslant 0}D_n'(t)\frac{x^n}{n!},\\
V=&\displaystyle
[4t-(1+3t)^2x+2t(3+2t+3t^2)x^2]
\sum_{n\geqslant 0}D_{n+1}(t)\frac{x^n}{n!},\\
VI=&\displaystyle
(1-t)^2\sum_{n\geqslant 0}D_n(t)\frac{x^n}{n!},
\end{array}
$$
so that
$$
\begin{array}{rl}
0=[x^n](I+\cdots+VI)=&\displaystyle
-\frac{4t^2(1-t)^2D_n''(t)}{n!}+
\frac{4t^2(1-t)^2(1+t)D_{n-1}''(t)}{(n-1)!}\\
&\displaystyle+\frac{4t(1-t)D_{n+1}'(t)}{n!}-
\frac{4t(1-t)(1+3t)D_n'(t)}{(n-1)!}\\
&\displaystyle+\frac{8t^2(1-t)(1+t)D_{n-1}'(t)}{(n-2)!}-
\frac{D_{n+2}(t)}{n!}+\frac{(1+5t)D_{n+1}(t)}{(n-1)!}\\
&-\displaystyle\frac{4t(1+2t)D_n(t)}{(n-2)!}+
\frac{4t^2(1+t)D_{n-1}(t)}{(n-3)!}\\
&-\displaystyle\frac{4t(1-t)^2D_n'(t)}{n!}+
\frac{2t(1-t)^2(3+t)D_{n-1}'(t)}{(n-1)!}\\
&+\displaystyle\frac{4tD_{n+1}(t)}{n!}-
\frac{(1+3t)^2D_n(t)}{(n-1)!}\\
&+\displaystyle\frac{2t(3+2t+3t^2)D_{n-1}(t)}{(n-2)!}+
\frac{(1-t)^2D_n(t)}{n!},
\end{array}
$$
where $[x^n](I+\cdots+VI)$ denotes the coefficient of $x^n$
in $I+\cdots+VI$.
Multiplying through by $n!$ and collecting terms,
the theorem follows.
\end{proof}

The recurrence relation (\ref{dnrec3}) involves $D_n(t)$ and
its derivatives only.
By equating coefficients of $t^k$,
we have the following recurrence relation for $D_{n,k}$.

\begin{Cor}\label{dnrec4}
$$
\begin{array}{rl}
D_{n+2,k}=&(n+4k)D_{n+1,k}\\
&+(5n-4k+8)D_{n+1,k-1}\\
&+[(1-n)-4(1+n)-4k(k-1)]D_{n,k}\\
&+[-2(1+n+2n^2)+8(1-n)(k-1)+8(k-1)(k-2)]D_{n,k-1}\\
&+[(1-n-8n^2)-4(1-3n)(k-2)-4(k-2)(k-3)]D_{n,k-2}\\
&+[6nk+4nk(k-1)]D_{n-1,k}\\
&+[6n(n-1)+2n(4n-9)(k-1)-4n(k-1)(k-2)]D_{n-1,k-1}\\
&+[4n(n-1)^2+2n(k-2)-4n(k-2)(k-3)]D_{n-1,k-2}\\
&+[2n(1-3n+2n^2)+2n(5-4n)(k-3)+4n(k-3)(k-4)]D_{n-1,k-3}.
\end{array}
$$
\end{Cor}

\bigskip
Theorem~\ref{dnrec5} (or equivalently Corollary~\ref{dnrec4})
constitutes the solution to a problem concerning the recurrence relations
for the Eulerian polynomials $D_n(t)$ (or the Eulerian number $D_{n,k}$),
posed by Brenti \cite{b94}.
Note that some of the coefficients of $D_{n,k}$ are not linear in $n$
as well as in $k$. A direct combinatorial proof of Corollary~\ref{dnrec4}
is desired, however.

\bigskip
A closely related problem, also posed by Brenti,
is whether the Eulerian polynomials $D_n(t)$ have real zeros only.
That $\mathfrak{S}_n(t)$ and $B_n(t)$ having only real zeros follow easily
from their recurrence relations.
Given the complicated nature of the recurrence relation for $D_n(t)$,
the corresponding assertion for $D_n(t)$ need not follow in a similar manner.

\bibliography{main}
\bibliographystyle{plain}

\end{document}